\documentclass[12pt]{article}
\usepackage{amssymb}
\def\hang{\hangindent\parindent}
\def\tex#1{\indent\llap{[#1]\enspace}\ignorespaces}
\def\re{\par\hang\tex}

\def\a{\alpha}
\def\b{\beta}

\def\D{\Delta}
\def\UU{{\cal U}}

\def\LL{{\cal L}}
\def\FF{{\cal F}}

\def\g{\gamma}
\def\G{\Gamma}

\def\e{\epsilon}

\def\ssc{\scriptscriptstyle}

\def\cl{\centerline}

\def\vs{\vspace*}

\def\ra{\rangle}

\def\la{\langle}
\def\ni{\noindent}
\def\ptl{\partial}

\def\Z{\mathbb{Z}{\ssc\,}}

\def\C{\mathbb{C}{\ssc\,}}
\def\F{\mathbb{F}{\ssc\,}}
\begin{document}

\cl{{\Large \bf
 Quantization of generalized Virasoro-like algebras }\footnote{Supported by NSF grant Y2006A17 of Shandong Province,
   NSF J06P52 of Shandong Provincial Education Department, China,
  NSF grant 10171064 of China,
 "One hundred Talents Program" from University of Science and Technology of China}}
\vs{6pt}

\cl{ Guang'ai Song$^*$, \ Yucai Su$^{\dag}$ and Yuezhu Wu$^{\ddag}$}
\cl{\small $^{*\,}$College of Mathematics and Information Science,
Shandong Institute of \vs{-4pt}Business }

\cl{ and Technology, Yantai, Shandong 264005, China}

\cl{$^{\dag\,}$Department of Mathematics, University of Science and
Technology of China\vs{-4pt},} \cl{ Hefei 230026, China}

\cl{$^{\ddag\,}$Department of Mathematics, Qufu Normal University,
Qufu, Shandong 273165,  China }

\cl{\small E-mail: gasong@ccec.edu.cn} \vs{6pt}

{\lineskip5pt\noindent{\small{\bf Abstract.} In a recent paper by
the authors, Lie bialgebras structures of generalized Virasoro-like
type were considered. In this paper, the explicit formula of the
quantization of generalized Virasoro-like algebras is presented.

  }
\vs{6pt}}

\cl{\bf\S1. \
Introduction}\setcounter{section}{1}\setcounter{equation}{0}

In Hopf algebras or quantum groups theory, there are two standard
methods to yield new bialgebras from old ones, one is twisting the
product by a 2-cocycle but keeping the coproduct unchanged, another
is twisting the coproduct by a Drinfel'd twist element but keeping
the product unchanged. Constructing quantization of Lie bialgebras
is an important method to produce new quantum groups (cf.~[V], [PO],
etc). In the paper [M1] (cf.~[M2, M3]), a class of infinite
dimensional Lie bialgebras containing Virasoro algebras was
presented. This type Lie bialgebras was classified in [NT], and the
quantization of this type Lie algebras was determined in [G]. In the
paper [SS], Lie bialgebras structures of generalized Witt type were
classified. The quantization of this type algebras was considered in
[HW]. In the paper [WSS], the structures of Lie bialgebras of
generalized Virasoro-like algebras were determined. In the preset
paper, we will consider the quantization of this type Lie algebras.
\vskip10pt

\cl{\bf\S2.
Preliminaries}\setcounter{section}{2}\setcounter{equation}{0}

\ni{\bf \S2.1 Generalized Virasoro-like Lie bialgebras.} Let $\C$
denote the complex field, and let $\G$ be any nondegenerate additive
subgroup of $\C^2$ (namely $\G$ contains a $\C$-basis of $\C^2$ ).

{\bf Definition 2.1} The Lie algebras $\LL (\G)$ with basis
$\{L_{\a} , \ptl_1 , \ptl_2 | \a \in \G \backslash \{ 0\} \}$ and
bracket:

\cl{$[L_{\a} , L_{\b} ] = (\a_1 \b_2 - \a_2 \b_1 )L_{\a + \b },
[\ptl_i , L_{\a} ] = \a_i L_{\a_i} ,$ for $\a , \b \in \G
\backslash \{ 0\}, i= 1, 2$.}

\noindent is called a generalized Virasoro-like algebras.

Remark: We use the convention that if an undefined notation appears
in an expression, we always treat it as zero; for instance, $L_{\a}
= 0,$ if $\a = 0.$ In particular, when $\G = \Z^2$ the derived
subalgebras $[\LL (\Z^2 ), \LL (\Z^2) ] ={\rm span}\{L_{\a} | \a \in
\Z^2 \backslash \{0\} \}$ is the Virasoro-like algebras (cf. [LT,
MJ, ZZ]). The Lie algebras $\LL (\G)$ is closely related to the Lie
algebras of Block type (cf. [DZ, X, Z]) and the Lie algebras of
Cartan type $S$ (cf. [SX, X1, Z2]).

The following theorem is the main result in [WSS].

{\bf Theorem 2.2} Every Lie bialgebra structure on the Lie algebras
$\LL (\G)$ is a coboundary triangular Lie bialgebra, namely, it is
given by $r$-matrix $r - r^{21}$ with $r = T \otimes L_{\a} , T \in
$ span$\{ \a_1 , \a_2\} , L_{\a} \in \LL (\G)$.\vskip7pt

\ni{\bf \S2.2 Drinfel'd twisting.} Let $A$ be a unital $R$-algebra
(where $R$ is a ring ), for any element $x \in A, a \in R,$ we set
\begin{equation}
x_a ^{\la n\ra } = (x + a)(x+ a+1)\cdots (x+a + n-1),
\end{equation}
\begin{equation}
x_a ^{[n]} = (x+ a)(x+a-1)\cdots (x+a-n+1),
\end{equation}
\noindent where $n \in \Z ,$ and we denote $x^{\la n\ra } = x_0
^{\la n\ra } ,
 x^{[n]} = x_0 ^{[n]}$.

 The following lemma belongs to [G] and [GZ].

 {\bf Lemma 2.3} Let $\F$ be a field with char$\,\F = 0$, and $x$ be
 any element of a unital $\F$-algebras $A$, for $a, d \in \F, $ and
 $m, n, r \in \Z$, one has:
 \begin{eqnarray}
 x_a ^{\la m + n\ra } &\!\!\!=\!\!\!& x_a ^{\la m\ra } x_{m+a} ^{\la n\ra },
\\
x_a ^{[m + n]} &\!\!\!=\!\!\!& x_a ^{[m]} x_{a-m} ^{[n]}, \\
 x_a ^{[m]} &\!\!\!=\!\!\!&
x_{a-m+1} ^{\la m\ra } , \\
 \mbox{$\!\!\!\!\!\!\!\!\sum\limits_{m+n = r}$} \frac{(-1)^n }{m!n!} x_a
^{[m]} x_d ^{\la n\ra } &\!\!\!=\!\!\!&\Big(\!\begin{array}{c} a- d \\
r
 \end{array}\!\Big) = \frac{(a-d)(a-d+1) \cdots (a-d-r+1)}{r!},
\\
 \mbox{$\!\!\!\!\!\!\!\!\sum\limits_{m+n=r}$} \frac{(-1)^m }{m!n!} x_a ^{[m]} x_{d-m} ^{[n]} &\!\!\!=\!\!\!&
 \Big(\!\begin{array}{c} a-d+r-1 \\ r \end{array}\! \Big) = \frac{(a-d)(a-d+1) \cdots (a-d
 +r-1)} {r!}.
 \end{eqnarray}
The following definition belongs to [D].

{\bf Definition 2.4} An element $\cal{F} \in H \otimes H$ is
called Drinfel'd twisting element, if it is invertible such that
\begin{eqnarray}
(\FF \otimes 1) (\D_0 \otimes Id ) (\FF ) &\!\!\!=\!\!\!& (1 \otimes
\FF) (1 \otimes \D_0 ) (\FF ) \\
(\e_0 \otimes Id ) (\FF ) &\!\!\!=\!\!\!& 1 \otimes 1 = (Id \otimes
\e_0 ) (\FF) .\end{eqnarray}

The following theorem is well known (cf.[D]), it can be found in
any book of hopf algebras.

 {\bf Theorem 2.5} Let ($H, \mu, \tau, \D_0, \e_0, S_0$) be a hopf
 algebra over commutative ring, $\cal{F}$ be a Drinfel'd element
 of $H \otimes H$, then

 (1) $u = \mu (Id \otimes S_0 ) (\cal{F} )$ is a invertible
 element of $H \otimes H$ with $u^{-1} = \mu (S_0 \otimes Id ) (\cal{F}) .$

 (2) The algebras $(H, \mu, \tau, \D, \e, S )$ is a new hopf
 algebra if we keep the counit undeformed and define
 $\D : H \rightarrow H \otimes H, \ S: H \rightarrow H $ by:

 \cl{$\D (h) = \FF \D_0 (h) \FF^{-1} , \ \ \ \ \  S(h) = u S_0 (h) u^{-1} . $}

Let $(\UU (\LL (\G)), \mu,  \tau, \D_0 , \e_0, S_0 )$ be the
standard hopf algebra, i.e. \begin{eqnarray*}
&\!\!\!\!\!\!\!\!\!\!\!\!\!\!\!\!&
 \D_0 (L_{\a} ) = L_{\a} \otimes 1
+ 1 \otimes L_{ \a} ,\ \    \D_0 (\ptl_{i}) = \ptl_{i} \otimes 1 + 1
\otimes \ptl_{i} ,\\&\!\!\!\!\!\!\!\!\!\!\!\!\!\!\!\!&
 S(L_{\a}) = -L_{\a}, \ \ S(\ptl_{i}) = -\ptl_{i} ,\ \ \ \
\e_0 (L_{\a} ) =0,\ \  \e_0 (\ptl_i) =0, \end{eqnarray*} for $\a \in
\F \backslash \{ 0 \}, i = 1, 2. $ The main result of this paper is
the following theorem.

{\bf Theorem 2.6} Let $\LL(\G)$ be the generalized Virasoro-like
algebra over $\F$ with char$\,\F = 0$, choose $T = a_1 \ptl_1 + a_2
\ptl_2 \in $ span$\{\ptl_1 , \ptl_2 \}$ and $L_{\a} \in \LL(\G)$
with $[T, L_{\a} ] = L_{\a} ,$ then there exists a noncommutative
and noncocommuative hopf algebra structere($\UU (\LL (\G))[[t]],
\mu, \tau, \D, S, \e $) on $\UU (\LL (\G))[[t]]$ over $\F[[t]]$,
such that $\UU (\LL (\G))[[t]]/t\UU (\LL (\G))[[t]] = \UU (\LL
(\G))$, which preserves the product and the counit of $\UU (\LL
(\G))[[t]]$, but the coproduct and  antipode are defined by
\begin{eqnarray}
\D(L_{\b}) &\!\!\!=\!\!\!& L_{\b} \otimes (1 - L_{\a} t )^b +
\mbox{$\sum\limits_{i=0} ^{\infty} $}(-1 )^{i} T^i \otimes (1-
L_{\a} t )^{-i} L_{\b + i \a} c_i t^i,
\\
\D (\ptl_i) &\!\!\!=\!\!\!& \ptl_j \otimes 1 + 1 \otimes \ptl_j +
\a_j T^{\la 1\ra
} \otimes (1- L_{\a} t)^{-1} L_{\a} t, J=1,2. \\
S (L_{\b} ) &\!\!\!=\!\!\!& - (1- L_{\a} t)^{-b}\mbox{$
\sum\limits_{i=0 } ^{\infty}$} L_{\b + i\a
} c_i T_1 ^{\la i\ra } t^i, \\
S(\ptl_i ) &\!\!\!=\!\!\!&\a_j T (1- L_{\a} t)^{-1}( L_{\a} t -
L_{\a} ^2 t^2 ) - \ptl_j
 \end{eqnarray}
\noindent where $b = a_1 \b_1 + a_2 + \b_2,\, c_i = \frac{(\a_1 \b_2
- \a_2 \b_1 )^{i} }{i!},\, c_0 =1,\, j=1, 2.\vspace*{5pt}$

 \cl{\bf\S3.
Proof of the main
result\vspace*{0pt}}\setcounter{section}{3}\setcounter{equation}{0}
We shall divide the proof of the Theorem 2.6 into several lemmas.

{\bf Lemma 3.1} Let $T, L_{\a} \in \LL(\G)$ with $$T = a_1 \ptl_1 +
a_2 \ptl_2 ,\ \ \ \ [T, L_{\a} ] = [a_1 \ptl_1 + a_2 \ptl_2 , L_{\a}
] = L_{\a},$$ for some $a_1,a_2\in\F$ and any $\a = (\a_1 ,
\a_2)\in\G$. For any $\b = (\b_1 , \b_2 ) \in \G$, denote $b = a_1
\b_1 + a_2 \b_2$. The following equations hold in $\UU (\LL (\G))$
for $a\in\F,\,m,k\in\Z_+$:
\begin{eqnarray}\!\!\!\!\!&\!\!\!\!\!&
L_{\b} T_a ^{[m]} = T_{a-b} ^{[m]} \b_{\b}, \\\!\!\!\!\!&\!\!\!\!\!&
L_{\b} T_a ^{\la m\ra } = T_{a-b} ^{\la m\ra } L_{\b},
\\\!\!\!\!\!&\!\!\!\!\!&
L_{\a} ^k T_a ^{[m]} = T_{a-k } ^{[m]} L_{\a} ^{k} ,
\\\!\!\!\!\!&\!\!\!\!\!& L_{\a} ^k T_a ^{\la m\ra } = T_{a-k } ^{\la
m\ra } L_{\a} ^{k} ,
\\\!\!\!\!\!&\!\!\!\!\!&
\ptl_j ^k T_{a} ^{[m]} = T_{a} ^{[m]} \ptl_j ^k , j= 1, 2,
\\\!\!\!\!\!&\!\!\!\!\!& \ptl_j ^k T_{a} ^{\la m\ra } = T_{a} ^{\la m\ra }
\ptl_j ^k , j= 1, 2, \\\!\!\!\!\!&\!\!\!\!\!& L_{\b} L_{\g} ^m =
\mbox{$\sum\limits_{i=0} ^{m }$} (-1)^{i} \Big(
\begin{array}{c} m \\ i \end{array} \Big)  (\g_1 \b_2 - \g_2 \b_1
)^i L_{\g} ^{m-i} L_{\b + i\g }, \\\!\!\!\!\!&\!\!\!\!\!& \ptl_j
L_{\g} ^m = m \g_j L_{\g,} ^{m} + L_{\g} ^m \ptl_j ,  j=1, 2.
\end{eqnarray}

{\it Proof.} Since $$[T, L_{\b} ] = [a_1 \ptl_1 + a_2 \ptl_2 ,
L_{\b} ] = (a_1 \b_1 + a_2 \b_2 ) L_{\b} = b L_{\b} = T L_{\b} -
L_{\b} T,$$ we have $L_{\b} T = (T -b ) L_{\b},$ it easy to see that
(2.1) is true for $m=1.$ Suppose that (3.1) is true for $m$, then
for $m+1$ we have

$L_{\b} T_a ^{[m+1]} = L_{\b} T_{a} ^{[m]} (T + a -m ) = T_{a-b}
^{[m]} L_{\b} (T + a -m) = T_{a -b} ^{[m]} (T + a - b -m ) L_{\b}
= T_{a-b} ^{[m+1]} .$

\noindent By induction on $m$, (3.1) holds. Similarly, we can obtain
(3.2)--(3.6). For (3.7), we have:
 $$\begin{array}{lll}L_{\b} L_{\g} ^m &\!\!\!=\!\!\!& \sum\limits_{i=0} ^{m} (-1)^{i }
\Big(\begin{array}{c} m \\ i
\end{array}\Big) L_{\g} ^{m-i} ({\rm ad} L_{\g})^{i} (L_{\b}) \\&\!\!\!=\!\!\!&
\sum\limits_{i=0} ^{m} (-1)^{i } \Big(\begin{array}{c} m \\ i
\end{array}\Big) L_{\g} ^{m-i} L_{\b + i\g} (\g_1 \b_2 - \g_2
\b_1 )^i .\end{array}$$ The proof of formula (3.8) is similar to
that of (3.7), using the following fact: $${\rm ad}(L_{\g} )^i
(\ptl_j )= \Big\{ \begin{array}{ll} - \g_j L_{\g} , \ \ & \ \ {\rm
if } \ \ i= 1,\\ 0, \ \ &\ \ {\rm if }\ \ i > 1.\end{array}
\eqno\Box$$

Now for $a \in \F$, set
$$\begin{array}{lll}
\FF_a = \sum\limits_{i = 0 } ^{\infty } \frac{(-1)^{i} }{i! } T_a
^{[i] } \otimes L_{\a} ^i t^i ,\\[9pt]
 F_a = \sum\limits_{i = 0} ^{\infty}
\frac{(-1)^{i} }{i! } T_a ^{\la i\ra  } \otimes L_{\a} t^i , \\[7pt]
u_a = \mu \cdot (S_0 \otimes Id ) (F_a ), v_a  = \mu \cdot (Id
\otimes S_0 ) (\FF_a ).\end{array}$$ Write $\FF = \FF_0 ,\, F = F_0
,\, u = u_0 , \,v = v_0 .$ Since $S_0 (T_a ^{\la i\ra  } ) = (-1 )^i
T_{-a } ^{[i] },\, S_0 (L_{\a } ^{i} ) = (-1)^{i} L_{\a} ,$ we have:
$$\begin{array}{lll}
u_a = \mu (S_0 \otimes Id ) (\sum\limits_{i=0} ^{\infty}
\frac{1}{i!} T_{a} ^{\la i\ra } \otimes L_{\a} ^i t^i ) =
\sum\limits_{i=0} ^{\infty} \frac{(-1)^i }{i!} T_{-a} ^{[i]} L_{\a}
^i t^i ,\\[9pt] v_a = \mu (Id \otimes S_0 ) (\sum\limits_{i=0} ^{\infty}
\frac{(-1)}{i!} T_{a} ^{[i] } \otimes L_{\a} ^i t^i ) =
\sum\limits_{i=0} ^{\infty} \frac{1}{i!} T_{a} ^{[i]} L_{\a} ^i
t^i.\end{array}$$

{\bf Lemma 3.2.} For $a, d \in \F, $ one has

\cl{$\FF_a F_d = 1 \otimes (1 - L_{\a} t )^{(a- d) } , v_a u_d =
(1- L_{\a} t)^{-(a + d)} $.}

Therefore the elements $\FF_a , F_a , u_a , v_a$ are invertible
elements with $\FF_a ^{-1} = F_{a} ,\, u_a ^{-1} = v_{-a}$.

{\it Proof.} Using the formula (2.6) we have:
$$\begin{array}{lll}\FF_a F_d &\!\!\!=\!\!\!& (\sum\limits_{i=0} ^{\infty} \frac{(-1)
^i}{i!} T_a ^{[i]} \otimes L_{\a}^i t^i ) \cdot (\sum\limits_{j=0}
^{\infty} \frac{1}{j!} T_d ^{\la j\ra } \otimes L_{\a} ^j t^j )\\[9pt]&\!\!\!=\!\!\!&\sum\limits_{i, j =0} ^{\infty} \frac{(-1)^{i}}{i!j!} T_{a} ^{[i]
} T_d ^{\la j\ra } \otimes L_{\a} ^i L_{\a} ^j t^{i+j}
=
\sum\limits_{m=0} ^{\infty} \big(\sum\limits_{i+j = m}
\frac{(-1)^j}{i!j!} T_a
^{[i]} T_d ^{\la j\ra } \big) \otimes L_{\a} ^m t^m \\[9pt]
&\!\!\!=\!\!\!& \sum\limits_{m=0} ^{\infty} (-1)^m
\Big(\begin{array}{c} a-d
\\ m
\end{array}\Big) L_{\a} ^{m} t^m
=
1 \otimes (1- L_{\a} t)^{a-d} ,
\\[13pt]
v_a u_d &\!\!\!=\!\!\!& (\sum\limits_{i=0} ^{\infty} \frac{1}{i!}
T_{a} ^{[i]} L_{\a} ^i t^i )\sum\limits_{j=0} ^{\infty} \frac{(-1)^j
}{j!} T_{-d} ^{[j]} L_{\a} ^j t^j
=
\sum\limits_{i, j = 0} ^{\infty} \frac{1}{i!} T_a ^{[i]} L_{\a}
^i
t^i \frac{(-1)^{j}}{j!} T_{-d} ^{[j]} L_{\a} ^j t^j \\
&\!\!\!=\!\!\!& \sum\limits_{m=0} ^{\infty}
\sum\limits_{i+j=m}\frac{(-1)^j }{i!j! } T_a ^{[i]} T_{-d-i}
^{[j]} L_{\a} ^{i+j} t^{i+j} \ \ \ \ \ ({\rm from} \ (3.3)) \\
&\!\!\!=\!\!\!& \sum\limits_{m=0} ^{\infty} \Big(\begin{array}{c}
a+d+m-1\\ m
\end{array}\Big) L_{\a} ^{m} t^m \ \ \ ({\rm from}\ (2.7)) \\
&\!\!\!=\!\!\!& (1-L_{\a} t)^{-(a+d)}\vspace*{-15pt}. \end{array}$$
 \hfill $\Box$

{\bf Lemma 3.3. } For any nonnegative integer $m$, and any $a \in
\F$ we have $$\mbox{$\D_0 T^{[m]} = \sum\limits_{i=0} ^{m}
\Big(\begin{array}{c} m \\ i
\end{array} \Big) T_{-a}^{[i]} \otimes T_a ^{[m-i]} .$}$$
In particular, we have $\D_0 T^{[m]} = \sum_{i=0} ^{m}
\big(\begin{array}{c}\! m \\[-8pt] i
\end{array}\! \big) T ^{[i]} \otimes T ^{[m-i]}$.

{\it Proof.} Since $\D_0 (T) = T \otimes 1 + 1 \otimes T, $ it is
easy to see that the result is true for $m=1.$ Suppose it is  true
for $m$, then for $m+1$, we have:
$$\begin{array}{lll}\D_0 (T^{[m+1]}) &\!\!\!=\!\!\!& \D_0 (T^{[m]}) \D_0 (T-m))\\
&\!\!\!=\!\!\!& (\sum\limits_{i=0} ^m  \Big( \begin{array}{c} m\\ i
\end{array} \Big) T_{-a} ^{[i]} \otimes T_a ^{[m-i]} ) ((T- a- m) \otimes 1 + 1 \otimes (T +a - m) + m (1 \otimes 1)
) \\ &\!\!\!=\!\!\!& (\sum\limits_{i = 1} ^{m-1} \Big(
\begin{array}{c} m\\ i
\end{array} \Big) T_{-a} ^{[i]} \otimes T_{a} ^{[m-i]} ) ((T -a- m) \otimes 1 + 1 \otimes (T +a- m) )  \\
&&+ m (\sum\limits_{i = 0} ^{m} \Big( \begin{array}{c} m\\ i
\end{array} \Big) T_{-a} ^{[i]} \otimes T_{a} ^{[m-i]}) + (1 \otimes T_a ^{[m+1]} + T_{-a}^{[m+1]} \otimes 1) \\
&& + (T -a - m) \otimes T_a ^{[m]} + T_{-a} ^{[m]} \otimes (T + a - m) \\
&\!\!\!=\!\!\!& 1 \otimes T_a ^{[m+1]} + T_{-a} ^{[m+1]} \otimes 1 +
m (\sum\limits_{i = 1} ^{m -1} \Big( \begin{array}{c} m\\ i
\end{array} \Big) T_{-a} ^{[i]} \otimes T_{a} ^{[m-i]})  \\
&&+ (T- a) \otimes T_a ^{[m]} + T_{-a} ^{[m]} \otimes (T +a) +
\sum\limits_{i=1} ^{m-1} \Big(\begin{array}{c} m\\ i \end{array}
\Big) T_{-a} ^{[i+1]} \otimes
T_a ^{[m-i]} \\
&&+ \sum\limits_{i=1} ^{m-1} (i-m)\Big( \begin{array}{c} m\\ i
\end{array} \Big) T_{-a} ^{[i]} \otimes T_a ^{[m-i]} + \sum\limits_{i=1} ^{m-1}
\Big( \begin{array}{c} m\\ i \end{array} \Big)T_{-a} ^{[i]} \otimes
T_a ^{[m-i +1]} \\
&&+ \sum\limits_{i=1} ^{m-1} (-i)\Big( \begin{array}{c} m\\ i
\end{array} \Big) T_{-a} ^{[i]} \otimes T_a ^{[m-i]}  \\
&\!\!\!=\!\!\!& 1 \otimes T_a ^{[m+1]} + T_{-a} ^{[m+1]} \otimes 1 +
\sum\limits_{i=1} ^{m} \Big( \Big( \begin{array}{c} m\\ i-1
\end{array} \Big) + \Big( \begin{array}{c} m\\ i
\end{array} \Big)   \Big) T_{-a} ^{[i]} \otimes T_a ^{[m+1-i]} \\
&\!\!\!=\!\!\!& \sum\limits_{i=0} ^{m+1} \Big( \begin{array}{c} m
+1\\ i
\end{array} \Big)T_{-a} ^{[i]} \otimes T_a ^{[m+1-i]} .
\end{array}$$
 From induction the result  hold on arbitrary $m$. \hfill $\Box$

{\bf Lemma 3.3.} $\FF = \sum_{i=0} ^{\infty} \frac{(-1)^i }{i!}
T^{[i]} \otimes L_{\a} ^i t^i $ is a Drinfel'd twist element of $\UU
(\LL(\G))[[t]]$, i.e.
$$\begin{array}{ll}
(\FF \otimes 1) (\D_0 \otimes Id ) (\FF ) = (1 \otimes \FF) (1
\otimes \D_0 ) (\FF ),\\[7pt]
(\e_0 \otimes Id ) (\FF ) = 1 \otimes 1 = (Id \otimes \e_0 ) (\FF)
.\end{array}$$

{\it Proof.} The second equation holds obviously, we just need to
prove the first one. Since $$\begin{array}{lll}(\FF \otimes 1) (\D_0
\otimes Id ) (\FF ) &\!\!\!=\!\!\!& (\sum\limits_{i=0} ^{\infty}
\frac{(-1)^i }{i!} T^{[i]} \otimes L_{\a} ^i t^i \otimes 1) (\D_0
\otimes Id ) (\sum\limits_{j=0} ^{\infty}
\frac{(-1)^j }{j!} T^{[j]} \otimes L_{\a} ^j t^j ) \\
&\!\!\!=\!\!\!&(\sum\limits_{i=0} ^{\infty} \frac{(-1)^i }{i!}
T^{[i]} \otimes L_{\a} ^i t^i \otimes 1) \times \\ &&
\times(\sum\limits_{j=0} ^{\infty} \frac{(-1)^j }{j!}
\sum\limits_{k=0} ^j \Big(\begin{array}{c} j\\ k
\end{array} \Big) T_{-i} ^{[k]} \otimes T_i ^{[j-k]} \otimes L_{\a}
^j t^j )\\ &\!\!\!=\!\!\!& \sum\limits_{i,j=0} ^{\infty}
\frac{(-1)^{i+j} }{i!j!} t^{i+j}\sum\limits_{k=0} ^j
\Big(\begin{array}{c} j\\ k
\end{array} \Big) T^{[i]} T_{-i}
^{[k]} \otimes L_{\a} ^i T_i ^{[j-k]} \otimes L_{\a} ^j \\
&\!\!\!=\!\!\!& \sum\limits_{i,j=0} ^{\infty}
\frac{(-1)^{i+j}}{i!j!} t^{i+j} \sum\limits_{k=0} ^j
\Big(\begin{array}{c} j\\ k \end{array} \Big) T^{[i+k]} \otimes
T^{[j-k]} L_{\a} ^i \otimes L_{\a} ^j ,
\end{array}$$
\noindent on the other hand, $$\begin{array}{lll}(1 \otimes \FF )(Id
\otimes \D_0 ) (\FF )&\!\!\!=\!\!\!& (\sum\limits_{r=0} ^{\infty}
\frac{(-1)^r
}{r!} t^r 1\otimes T^{[r] } \otimes L_{\a} ^r  ) \times \\
&&\times(\sum\limits_{s=0 } ^{ \infty} \frac{(-1)^s }{s!} t^s T^{[s]
} \otimes \sum\limits_{q=0} ^s \Big(
\begin{array}{c} s\\ q
\end{array} \Big) L_{\a} ^q \otimes L_{\a} ^{s-q} ) \\
&\!\!\!=\!\!\!& \sum\limits_{r,s=0} ^{\infty} \frac{(-1)^{r+s}
}{r!s!} t^{r+s} \sum\limits_{q=0} ^s \Big(
\begin{array}{c} s\\ q
\end{array} \Big) T^{[s] } \otimes T^{[r] } L_{\a} ^q \otimes L_{\a}
^{r+s-q} ,
\end{array}$$
\noindent it sufficient to show for a fixed $m$ that:
$$\mbox{$ \sum\limits_{i+j=m}  \frac{1}{i!j!} t^{i+j} \sum\limits_{k=0} ^j
\Big(\begin{array}{c} j\\ k \end{array} \Big) T^{[i+k]} \otimes
T^{[j-k]} L_{\a} ^i \otimes L_{\a} ^j
=\sum\limits_{r+s=m}  \frac{1}{r!s!} t^{r+s}
\sum\limits_{q=0} ^s \Big(
\begin{array}{c} s\\ q
\end{array} \Big) T^{[s] } \otimes T^{[r] } L_{\a} ^q \otimes L_{\a}
^{r+s-q} .$}$$
 Now fix $r,\, 0 \leq i \leq s,$ and set $i=q, i+k = s,$ then we have
$$j-k = j-(s-i) = i + j - s =m-s =r. $$ We see that the
coefficients of $T^{[s]} \otimes T^{[r]} L_{\a} ^q \otimes L_{\a}
^{m-q}$ in both sides are equal. So the result holds. \hfill $\Box$

{\bf Lemma 3.4.} For $a \in \F, \b \in \G,$ we have
\begin{eqnarray} (L_{\b} \otimes 1 )F_a  &\!\!\!=\!\!\!& F_{a-b} (L_{\b}  \otimes
1), \\
 (1 \otimes L_{\b} ) F_a &\!\!\!=\!\!\!& \mbox{$\sum\limits_{l=0}
^{\infty}$} (-1)^l F_{a+l} (T_a ^{\la l\ra } \otimes c_l L_{\b+ l\a
}
t^l), \\
L_{\b} u_a &\!\!\!=\!\!\!& u_{a+b} \mbox{$\sum\limits_{l=0}
^{\infty}$} L_{\b
+ la} c_l T_{-a+l} ^{\la l\ra } t^l , \\
(\ptl_j  \otimes 1) F_a &\!\!\!=\!\!\!&F_a (\ptl_j \otimes 1), \\
(1 \otimes \ptl_j) F_a &\!\!\!=\!\!\!& F_{a+1} T_a ^{\la 1\ra }
\otimes \a_j L_{\a} t) + F_a (1 \otimes \ptl_j), \\
\ptl_j u_a &\!\!\!=\!\!\!&- \a_j T_{-a} ^{[1]} u_{a+1} L_{\a} t +
u_a \ptl_j , \\
L_{\a} u_a &\!\!\!=\!\!\!& u_{a+1} L_{\a} , \\
v T_{-a} ^{[1]} &\!\!\!=\!\!\!& T_{-a} ^{[1]} v_{-a} - T_{a} ^{[1]}
v_a L_{\a} t.
\end{eqnarray}
\noindent where $c_l = \frac{1}{l!} (\a_1 \b_2 - \a_2 \b_1 )^l , j
=1,2 .$

{\it  Proof.} From (3.2) we have:
$$\begin{array}{lll}(L_{\b} \otimes 1) F_a &\!\!\!=\!\!\!& (L_{\b} \otimes 1)
\sum\limits_{i=0} ^{\infty} \frac{1}{i!} T_a ^{\la i\ra } \otimes
L_{\a} ^i t^i \\
&\!\!\!=\!\!\!& \sum\limits_{i=0} ^{\infty} \frac{1}{i!} L_{\b} T_a
^{\la i\ra } \otimes
L_{\a} ^i t^i \\[6pt]
&\!\!\!=\!\!\!& \sum\limits_{i=0} ^{\infty} \frac{1}{i!} T_{a-b}
^{\la i\ra } L_{\b} \otimes L_{\a} ^i t^i
=
F_{a-b} (L_{\b} \otimes 1),
\end{array}$$
\noindent where $b = \a_1 \b_2 + \a_2 \b_1 .$ This proves (3.9).
 For (3.10), using (3.7) we have:
$$\begin{array}{lll} (1 \otimes L_{\b} ) F_a &\!\!\!=\!\!\!& (1 \otimes L_{\b})
 \sum\limits_{i=0} ^{\infty} \frac{1}{i!} T_a ^{\la i\ra } \otimes L_{\a} ^i t^i
 =
 \sum\limits_{i=0} ^{\infty} \frac{1}{i!} T_a ^{\la i\ra } \otimes L_{\b} L_{\a} ^i t^i \\
&\!\!\!=\!\!\!& \sum\limits_{i=0} ^{\infty} \frac{1}{i!} T_a ^{\la
i\ra } \otimes \sum\limits_{l=0} ^{i} (-1)^l \Big(\begin{array}{c}
i\\ l
\end{array} \Big) L_{\a} ^{i-l} L_{\b + l\a } (\a_1 \b_2 - \a_2
\b_1)^l t^i
\\[9pt] &\!\!\!=\!\!\!& \sum\limits_{i=0} ^{\infty}(\sum\limits_{l=0} ^i (-1)^l \frac{1}{(i-l)!l!} T_a ^{\la i\ra }
 \otimes L_{\a} ^{i-l} L_{\b + l\a }(\a_1 \b_2 - \a_2 \b_1)^l )
 t^i \\[9pt] &\!\!\!=\!\!\!& \sum\limits_{i=0} ^{\infty} \sum\limits_{l=0} ^{
 \infty} (-1)^l \frac{1}{i! l!} T_a ^{\la i+l\ra } \otimes L_{\a} ^i L_{\b + l\a
 } (\a_1 \b_2 - \a_2 \b_1 )^l t^{i+l} \\[9pt]
 &\!\!\!=\!\!\!& \sum\limits_{l=0} ^{\infty} (-1)^l \sum\limits_{i=0} ^{\infty} (\frac{1}{i!} T_{a+l} ^{\la i\ra } \otimes L_{\a} ^i t^i )
 \frac{1}{l!} T_{a} ^{\la l\ra } \otimes L_{\b + l\a } t^l (\a_1 \b_2 - \a_2 \b_1
 )^l \\[9pt] &\!\!\!=\!\!\!& \sum\limits_{l=0} ^{\infty } (-1)^l F_{a+l} (T_a ^{\la l\ra } \otimes L_{\b + l\a } c_l
 t^l). \end{array}$$
\noindent So we have (3.10).
 Now we prove (3.11):

$\begin{array}{lll} L_{\b} u_a &\!\!\!=\!\!\!& L_{\b}
\sum\limits_{r=0} ^{\infty}
  \frac{(-1)^r}{r!} T_{-a} ^{[r]} L_{\a} ^r t^r
  =
  \sum\limits_{r=0} ^{\infty}
  \frac{(-1)^r}{r!} L_{\b} T_{-a} ^{[r]} L_{\a} ^r t^r
  =
  \sum\limits_{r=0} ^{\infty}
  \frac{(-1)^r}{r!}  T_{-a -b} ^{[r]} L_{\b} L_{\a} ^r t^r \\
&\!\!\!=\!\!\!&\sum\limits_{r=0} ^{\infty}
  \frac{(-1)^r}{r!}  T_{-a -b} ^{[r]} \sum\limits_{l=0} ^r (-1)^l \Big(\begin{array}{c} r\\ l \end{array} \Big)
  L_{\a} ^{r-l} (\a_1 \b_2 - \a_2 \b_1 )^l L_{\b + l\a } t^r \\
  &\!\!\!=\!\!\!& \sum\limits_{r=0} ^{\infty} \frac{(-1)^r }{r!} T_{-a-b} ^{[r]}
  \sum\limits_{l=0} ^{r} (-1)^l \frac{r!}{(r-l)! \ l! } L_{\a} ^{r-l} (\a_1 \b_2 - \a_2
  \b_1)^l L_{\b + l\a } t^r
\end{array}$

$\begin{array}{lll}\phantom{L_{\b} u_a }
  &\!\!\!=\!\!\!& \sum\limits_{r, l =0} ^{\infty} \frac{(-1)^r }{r! \ l!} T_{-a-b}
  ^{[r+l]} L_{\a} ^r (\a_1 \b_2 - \a_2 \b_1)^l L_{\b + l\a } t^{r+l}
  \\[9pt]&\!\!\!=\!\!\!& \sum\limits_{r, l =0} ^{\infty} \frac{(-1)^r }{r! \ l!} T_{-a-b}
  ^{[r]} T_{-a-b-r} ^{[l]} L_{\a} ^r (\a_1 \b_2 - \a_2 \b_1)^l L_{\b + l\a }
  t^{r+l} \\[9pt]
  &\!\!\!=\!\!\!&  \sum\limits_{l =0} ^{\infty} \sum\limits_{r=0} ^{\infty}( \frac{(-1)^r }{r! } T_{-a-b}
  ^{[r]} L_{\a} ^r t^r )T_{-a-b} ^{[l]} L_{\b +l\a }\frac{(\a_1 \b_2 - \a_2 \b_1)^l }{l!}
  t^{l} \\[9pt]
  &\!\!\!=\!\!\!& u_{a+b} \sum\limits_{l=0} ^{\infty} T_{-a-b} ^{[l]} L_{\b +l\a }
  c_l t^l
=
u_{a+b} \sum\limits_{l=0} ^{\infty} L_{\b + l\a } T_{-a + l} ^{[l]}
  c_l t^l
=
u_{a+b} \sum\limits_{l=0} ^{\infty} L_{\b + l\a } T_{1-a} ^{\la l\ra }
  c_l t^l .
 \end{array}$

 \noindent This proves (3.11). For (3.12) we have:
$$\begin{array}{lll} (\ptl_j \otimes 1 )F_a &\!\!\!=\!\!\!&( \ptl_j \otimes
1)\sum\limits_{i=0} ^{\infty} \frac{1}{i!}
  T_a ^{\la i\ra } \otimes L_{\a} ^i t^i
=
\sum\limits_{i=0} ^{\infty} \frac{1}{i!}
  \ptl_j T_a ^{\la i\ra } \otimes L_{\a} ^i t^i \\
  &\!\!\!=\!\!\!&\sum\limits_{i=0} ^{\infty} \frac{1}{i!}
   T_a ^{\la i\ra } \ptl_j \otimes L_{\a} ^i t^i
=
(\sum\limits_{i=0} ^{\infty} \frac{1}{i!}
   T_a ^{\la i\ra } \otimes L_{\a} ^i t^i ) (\ptl_j \otimes 1)
=
F_a (\ptl_j \otimes 1).
\end{array}$$
Using (2.3) and (3.8) we have:
$$\begin{array}{lll}(1 \otimes \ptl_j ) F_a &\!\!\!=\!\!\!& (1 \otimes \ptl_j)
\sum\limits_{i=0} ^{\infty} \frac{1}{i!} T_a ^{\la i\ra } \otimes
L_{\a} ^i t^i
=
\sum\limits_{i=0} ^{\infty} \frac{1}{i!} T_a ^{\la i\ra }
\otimes\ptl_j L_{\a} ^i t^i\\ &\!\!\!=\!\!\!&\sum\limits_{i=0}
^{\infty} \frac{1}{i!} T_a ^{\la i\ra } \otimes (i \a_j L_{\a} ^i +
L_{\a} ^i \ptl_j) t^i
\\[9pt]&\!\!\!=\!\!\!&
\sum\limits_{i=0} ^{\infty} \frac{1}{(i-1)!} T_a ^{\la i\ra }
\otimes \a_j L_{\a} ^i t^i + \sum\limits_{i=0} ^{\infty}
\frac{1}{i!} T_a ^{\la i\ra } \otimes  L_{\a} ^i \ptl_j t^i\\[9pt]
 &\!\!\!=\!\!\!& \sum\limits_{i=0} ^{\infty}\frac{1}{(i -1)!} T_a ^{\la 1\ra } T_{a+1}
^{\la i-1\ra }\otimes \a_j L_{\a} ^i t^i + F_a (1 \otimes \ptl_j)
\\&\!\!\!=\!\!\!&
F_{a+1} (T_a ^{\la 1\ra } \otimes \a_j L_{\a} t) + F_a (1
\otimes \ptl_j).
\end{array}$$
So (3.13 ) holds. For (3.14) we have:
$$\begin{array}{lll}\ptl_j u_a &\!\!\!=\!\!\!& \ptl_j \sum\limits_{r=0} ^{\infty}
  \frac{(-1)^r}{r!} T_{-a} ^{[r]} L_{\a} ^r t^r
=
\sum\limits_{r=0} ^{\infty}
  \frac{(-1)^r}{r!} T_{-a} ^{[r]}\ptl_j L_{\a} ^r t^r \\[9pt]
&\!\!\!=\!\!\!&\sum\limits_{r=0} ^{\infty}
  \frac{(-1)^r}{r!} T_{-a} ^{[r]} (r \a_j L_{a} ^r  + L_{\a} ^r \ptl_j) t^r
  \\[9pt]
&\!\!\!=\!\!\!&\sum\limits_{r=0} ^{\infty} \a_j T_{-a} ^{[1]}
  \frac{(-1)^r}{(r-1)!} T_{-a-1} ^{[r-1]}   L_{a}
  ^r t^r +\sum\limits_{r=0} ^{\infty}
  \frac{(-1)^r}{r!} T_{-a} ^{[r]} L_{\a} ^r t^r \ptl_j \\
&\!\!\!=\!\!\!&-\a_j T_{-a} ^{[1]} u_{a+1} L_{\a} t + u_a \ptl_j .
\end{array}$$
Now we prove (3.15): $$\begin{array}{lll} L_{\a} u_a
=
L_{\a}
\sum\limits_{i=0} ^{\infty}
\frac{(-1)^i }{i!} T_{-a} ^{[i]} L_{\a} ^i t^i
=
\sum\limits_{i=0} ^{\infty} \frac{(-1)^i }{i!} T_{-a-1} ^{[i]}
L_{\a}
^{i +1} t^i
=
u_{a+1} L_{\a} . \end{array}$$
 For the last equation we have:
$$\begin{array}{lll} v T_{-a} ^{[1]} &\!\!\!=\!\!\!& \sum\limits_{i=0} ^{ \infty}
\frac{1}{i!} T_a ^{[i]} L_{\a} ^i t^i T_{-a} ^{[1]}
=
\sum\limits_{i=0} ^{ \infty} \frac{1}{i!} T_a ^{[i]} (T -a -i)
L_{\a} ^i t^i \\
&\!\!\!=\!\!\!& T_{-a} ^{[1]} v_a - \sum\limits_{i=0} ^{ \infty}
\frac{1}{(i-1)!}
(T+a ) T_{a-1} ^{[i-1]} L_{\a} ^i t^i \\
&\!\!\!=\!\!\!& T_{-a} ^{[1]} v_a  - T_{a} ^{[1]} v_{a-1} L_{\a} t .
\end{array}$$
This complete the proof of the lemma.\hfill$\Box$

{\it Proof of Theorem 2.6.} For arbitrary elements $L_{\b} , \ptl_j
\in \LL(\G), j=1, 2 $, we have:

\hspace*{15pt}$\begin{array}{lll}\D(L_{\b} ) &\!\!\!=\!\!\!& \FF
\D_0 (L_{\b}) \FF^{-1}
=
\FF (L_{\b} \otimes 1) \FF^{-1} + \FF (1 \otimes
L_{\b})
\FF^{-1}
\\&\!\!\!=\!\!\!&
\FF (L_{\b} \otimes 1) F + \FF (1 \otimes L_{\b})
F \\
&\!\!\!=\!\!\!& \FF F_{-b} (L_{\b} \otimes 1 ) + \FF
\sum\limits_{l=0} ^{\infty}
(-1)^l F_l (T^{\la l\ra } \otimes L_{\b + l\a} c_l t^l) \\
&\!\!\!=\!\!\!& (1 \otimes (1-L_{\a} t)^b ) (L_{\b} \otimes 1) \\
&&+\sum\limits_{l=0} ^{\infty} (-1)^l (1 \otimes (1- L_{\a} t)^{-l}
) \otimes (
T^{\la l\ra } \otimes L_{\b + l\a} c_l t^l) \\
&\!\!\!=\!\!\!& L_{\b} \otimes (1-L_{\a} t)^b + \sum\limits_{l=0}
^{\infty} (-1)^l T^{\la l\ra } \otimes (1-L_{\b} t)^{-l} L_{\b+l\a }
c_l t^l ,
\end{array}$

\hspace*{15pt}$\begin{array}{lll} \D (\ptl_j ) &\!\!\!=\!\!\!& \FF
\D(\ptl_j) \FF^{-1}
=
\FF (\ptl_j \otimes 1 + 1 \otimes \ptl_j ) F \\
&\!\!\!=\!\!\!& \FF (\ptl_j \otimes 1 ) F + \FF (1 \otimes \ptl_j ) F \\
&\!\!\!=\!\!\!& \FF F (\ptl_j \otimes 1 ) + \FF ( F_{1} (T ^{\la
1\ra } \otimes \a_j
L_{\a} t)+ F (1 \otimes \ptl_j) ) \\
&\!\!\!=\!\!\!& \ptl_j \otimes 1 + 1\otimes \ptl_j + 1 \otimes (1 -
L_{\a}
t)^{-1} (T ^{\la 1\ra } \otimes \a_j L_{\a} t)\\
&\!\!\!=\!\!\!&\ptl_j \otimes 1 + 1\otimes \ptl_j + \a_j T ^{\la
1\ra } \otimes (1 - L_{\a} t)^{-1}  L_{\a} t, j=1, 2,
\end{array}$

\hspace*{15pt}$\begin{array}{lll}S(L_{\b}) &\!\!\!=\!\!\!& u^{-1} S_{0} (L_{\b}) u
=
- v L_{\b} u
=
-v u_b (\sum\limits_{l=0} ^{\infty}c_l L_{\b + l\a } T_1 ^{\la l\ra
} t^l) \\ &\!\!\!=\!\!\!& -(1- L_{\a} t)^{-b}(\sum\limits_{l=0}
^{\infty}c_l L_{\b + l\a } T_1 ^{\la l\ra } t^l) ,
\end{array}$

\hspace*{15pt}$\begin{array}{lll} S(\ptl_j) &\!\!\!=\!\!\!& u^{-1} S_0 (\ptl) u
=
- v
\ptl_j u \\ &\!\!\!=\!\!\!& -v (- \a_j T ^{[1]} u_1 L_{\a} t + u \ptl_j) \\
&\!\!\!=\!\!\!& \a_j (T v -T v L_{\a } t) u_1 L_{\a} t -\ptl_j \\
&\!\!\!=\!\!\!& \a_j T v u_1 L_{\a} t - \ptl_j T v u_2 L_{\a} ^2 t^2  - \ptl_j\\
&\!\!\!=\!\!\!& \a_j T (1- L_{\a} t)^{-1} L_{\a} t - \a_j T
(1-L_{\a} t)^{-1}
L_{\a} ^2 t^2 - \ptl_j \\
&\!\!\!=\!\!\!& \a_j T (1- L_{\a} t)^{-1} (L_{\a} t - L_{\a} ^2 t^2)
- \ptl_j.
 \ \ \ \ \ \ \ \ \ \ \ \ \ \ \ \ \ \ \ \ \ \ \ \ \ \ \ \ \ \ \ \ \ \
\end{array} $\hfill $\Box$

\vskip9pt\ni{\bf References}\vskip5pt\small
\parindent=8ex\parskip=2pt\baselineskip=2pt

\re{D1} V. Drinfel'd, Constant quasiclassical solutions of the
Yang-Baxter quantum equation, {\it Soviet Math. Dokl.} {\bf28}(3)
(1983), 667--671.

\re{D2} V. Drinfel'd Quantum groups, Proceedings ICM (Berkeley
1986), Providence, Amer Math Soc, 1987, 789--820.

\re{F} R. Farnsteiner, Derivations and central extensions of
finitely generalized Lie algebras, {\it J. Algebra} {\bf 118}
(1988), 33--45.

\re{GZ} A. Giaquinto, J. Zhang, Bialgebra action, twists and
universal deformation formulas, {\it J. Pure Appl. Algebra}
{\bf128}(2) (1998), 133--151.

\re{G} C. Grunspan, Quantizations of the Witt algebra and of simple
Lie algebras in characteristic $p$, {\it J. Algebra} {\bf280}
(2004), 145--161.

\re{HW} N. Hu, X. Wang, Quantizations of generalized-Witt algebra
and of Jacobson-Witt algebra in modular case, arXiv:math.QA/0602281.

\re{M1} W. Michaelis, A class of infinite-dimensional Lie bialgebras
containing the Virasoro algebras, {\it Adv. Math.} {\bf107} (1994),
365--392.

\re{M2} W. Michaelis, Lie coalgebras, {\it Adv. Math.} {\bf38}
(1980), 1--54.

\re{M3} W. Michaelis, The dual Poincare-Birkhoff-Witt theorem, {\it
Adv. Math.} {\bf57} (1985), 93--162.

\re{NT} S.-H. Ng, Earl J. Taft, Classification of the Lie bialgebra
structures on the Witt and Virasoro algebras, {\it J. P. App.
algebra.} {\bf151} (2000), 67--88.

\re{N} W.D. Nichols, The structure of the dual Lie coalgebra of the
Witt algebra, {\it J. P. APP. Alg.} {\bf 68} (1990), 395--364W

\re{P}  D. Passman, New Simple Infinite Dimensional Lie Algebras,
{\it J. Algebra} {\bf 206} (1998), 682--692.

\re{SS} G. Song, Y. Su, Lie bialgebras of generalized Witt type,
{\it Science in China, Series A--Mathematics} {\bf 49}(4) (2006),
533--544.

\re{T} Earl J. Taft, Witt and Virasoro algebras as Lie bialgebras,
{\it J. P. App. algebra. } {\bf87} (1993), 301--312.

\re{WSS} Y. Wu, G. song, Y. Su, Lie bialgebras of generalized
Virasoro-like type, {\it Acta Mathematica Sinica, English Series}
{\bf22}(6) (2006), 1915--1922.
\end{document}